\documentstyle{amsart}
\font\tenmsb=msbm10
\font\sevenmsb=msbm7
\font\fivemsb=msbm5
\font\largemsb=msbm10 at11pt
\font\largesevenmsb=msbm7 at 7.7pt
\font\largefivemsb=msbm5 at 5.5pt
\font\smallmsb=msbm7 at7.7pt
\newfam\msbfam
\textfont\msbfam=\tenmsb
\scriptfont\msbfam=\sevenmsb
\scriptscriptfont\msbfam=\fivemsb
\newfam\largemsbfam
\textfont\largemsbfam=\largemsb
\scriptfont\largemsbfam=\largesevenmsb
\scriptscriptfont\largemsbfam=\largefivemsb
\newfam\smallmsbfam
\textfont\smallmsbfam=\smallmsb
\def\Bbb#1{\ifdim1em>10.1pt{\fam\largemsbfam\relax#1\kern.9pt}
   \else\ifdim1em<9.9pt{\fam\smallmsbfam\relax#1\kern.5pt}
   \else\fam\msbfam\relax#1\kern.8pt\fi\fi}
\def\qed{\hfill\vbox{\hrule height.09ex
   \hbox{\vrule width.09ex height1.7ex depth.3ex \kern2ex
   \vrule width.09ex height1.7ex depth.3ex}\hrule height.09ex}\bigskip}
\def\BOX{\hspace{0.015in}\hfill\vbox{\hrule height.09ex
   \hbox{\vrule width.09ex height1.4ex depth.3ex \kern1.8ex
   \vrule width.07ex height1.4ex depth.3ex}\hrule height.07ex}\hspace{0.02in}}

\def\R{{\Bbb R}}
\def\C{{\Bbb C}}

\def\D{{\Bbb D}}
\newcommand{\dis}{\displaystyle}
\def\mapright#1{\smash{
	\mathop{\stackrel{#1}{\rightarrow}}}}
\def\mapdown#1{\downarrow
	\rlap{$\vcenter{\hbox{$\scriptstyle#1$}}$}}

\def\DATE{
}

\newtheorem{theorem}{Theorem}
\newtheorem{proposition}{Proposition}[section]
\newtheorem{lemma}[proposition]{Lemma}

\begin{document}

\large

\title[On a boundary value problem]
{On a boundary value problem in subsonic aeroelasticity and
the cofinite Hilbert transform.}
\author[Peter L. Polyakov]{Peter L. Polyakov}
\address{Department of Mathematics, University of Wyoming, Laramie, WY 82071, USA}
\email{ polyakov@@uwyo.edu}
\date{\today}

\begin{abstract}
We study a boundary value problem in subsonic aeroelasticity and introduce the
{\it cofinite Hilbert transform} as a tool in solving an auxiliary linear integral
equation on the complement of a finite interval on the real line $\R$.
\end{abstract}
\maketitle

\centerline{\DATE}

\bigskip

\section{Introduction.}\label{Introduction}
\indent
We consider the linearized subsonic inviscid compressible flow equation in 2D
(\cite{BAH}, \cite{Ba2})
\begin{equation}\label{LinearEquation}
a^2_{\infty}\left(1-M^2\right)\frac{\partial^2\phi}{\partial x^2}
+a^2_{\infty}\frac{\partial^2\phi}{\partial z^2}
=\frac{\partial^2\phi}{\partial t^2}
+2Ma_{\infty}\frac{\partial^2\phi}{\partial t\partial x},
\end{equation}
where $a_{\infty}$ is the speed of sound,
${\dis M=\frac{U}{a_{\infty}}<1}$ - the Mach number, $U$ - free stream velocity,
$\phi(x,z,t)$ - small disturbance velocity potential, considered on
$$\R^2_{+}\times\overline{\R_{+}}=\left\{(x,z,t):-\infty<x<\infty,\ 0<z<\infty,\
0\leq t<\infty\right\},$$
with boundary conditions:
\begin{itemize}
\item
flow tangency condition
\begin{equation}\label{Flow_Tangency}
\frac{\partial\phi}{\partial z}(x,0,t)=w_a(x,t),\ |x|<b,
\end{equation}
where $b$ is the "half-chord", and $w_a$ is the given normal velocity of the wing,
without loss of generality we will assume in what follows that $b=1$,
\item
{\it "strong Kutta-Joukowski condition"} for the acceleration potential
$$\psi(x,z,t):=\frac{\partial\phi}{\partial t}+U\frac{\partial\phi}{\partial x},$$
\begin{equation}\label{Kutta}
\psi(x,0,t)=0\ \mbox{for}\ 1<|x|<A\ \mbox{for some}\ A>1,
\end{equation}
\item
far field condition
$$\phi(x,z,t)\to 0,\ \mbox{as}\ |x|\to\infty,\ \mbox{or}\ z\to\infty.$$
\end{itemize}
\indent
Boundary condition (\ref{Kutta}), though being motivated by one of the "auxiliary boundary
conditions" from (\cite{BAH}, p. 319), is weaker, because it requires
that $\psi(x,0,t)=0$ not on the whole $\R\setminus[-1,1]$, but only on finite intervals
adjacent to the interval $[-1,1]$. On the other hand this change in the boundary condition
allows application of some new mathematical tools different from tools in \cite{BAH} and
\cite{Ba2}.\\

\indent
In order to formulate our main result we introduce the following notations. We
denote by $\widehat{w_a}$ the Laplace transform of the function $w_a$ with respect
to time variable
$$\widehat{w_a}(x,z,\lambda)=\int_0^{\infty}e^{-\lambda t}w_a(x,z,t)dt$$
for $\mbox{Re}\lambda>\sigma_a>0.$ We also denote
$$\begin{array}{ll}
{\dis r(\lambda)=\frac{\lambda M}{U\sqrt{1-M^2}},}
\vspace{0.1in}\\
{\dis d(\lambda)=\frac{\lambda M^2}{U(1-M^2)}.}
\end{array}$$
\indent
In sections \ref{Solvability} and \ref{GLambdaResolvent} we construct a function
${\cal D}_{N}(\lambda)$ (equation (\ref{D_N})), analytic in the half-plane
$\mbox{Re}\lambda>\sigma_a>0$, and depending only on the function
$K_0$ - the modified Bessel function of the third kind.

\indent
The following theorem represents the main result of the paper.

\begin{theorem}\label{Main}\ Let function ${\cal D}_{N}(\lambda)$ from equation (\ref{D_N}),
mentioned above, have no zeros in the strip $\left\{\mbox{Re}\lambda\in [\sigma_1,\sigma_2]\right\}$,
where $\sigma_1>\sigma_a$. Let $I(1)=[-1,1]$, and let
$w_a(\cdot,t) \in L^2\left(I(1)\right)$ be such that for some $\epsilon>0$
\begin{equation}\label{wCondition}
{\dis \left\|\widehat{w}_a(\cdot,\sigma+i\eta)\right\|_{L^2\left(I(1)\right)}
<\exp\left\{-e^{|\eta|}\cdot(1+|\eta|)^{2+\epsilon}\right\}\ \mbox{for}\
\sigma\in[\sigma_1,\sigma_2] }
\end{equation}
\indent
Then equation (\ref{LinearEquation}) has a solution of the form
\begin{equation}\label{Solution}
\begin{array}{ll}
{\dis \phi(x,z,t)=-\frac{1}{2\pi\sqrt{1-M^2}}
\int_{\sigma^{\prime}-i\infty}^{\sigma^{\prime}+i\infty}e^{d(\sigma^{\prime}+i\eta)x} }
\vspace{0.1in}\\
{\dis \times\left[\int_{-\infty}^\infty K_0\left(r(\sigma^{\prime}+i\eta)
\left(\frac{(x-y)^2}{1-M^2}+z^2\right)^{\frac{1}{2}}\right)
h_a(y,\sigma^{\prime}+i\eta)dy\right]e^{(\sigma^{\prime}+i\eta)t}d\eta.}
\end{array}
\end{equation}
This solution is independent of $\sigma^{\prime}\in [\sigma_1,\sigma_2]$,
satisfies boundary conditions above, and function $h_a$ satisfies the estimate
$$\int_{-\infty}^\infty(1+|x|)^{p-2}\left|h_a(x,\sigma^{\prime}+i\eta)\right|^pdx
<\frac{C}{\left(1+|\eta|\right)^m}$$
for arbitrary $m>0$, ${\dis p<\frac{4}{3} }$, and $C>0$ independent of $\lambda$.
\end{theorem}
\indent
The author thanks A.V. Balakrishnan for suggesting and explaining the problem considered
here, and for his hospitality during author's visits to UCLA.

\section{"General" solution.}
\label{General}

\indent
We are seeking a solution of the equation (\ref{LinearEquation}) of the form
\begin{equation}\label{SolutionForm}
\phi(x,z,t)=\int_{\sigma-i\infty}^{\sigma+i\infty}
\xi(x,z,\lambda)e^{(\sigma+i\eta)t}d\eta,
\end{equation}
where $\lambda=\sigma+i\eta$, $\sigma>\sigma_a$ and $\xi(x,z,\lambda)\in L^1\eta(\R)$.
Then, substituting the expression above into equation (\ref{LinearEquation}),
we obtain the following auxiliary equation for $\xi$
\begin{equation}\label{LaplaceEquation}
a^2_{\infty}\left(1-M^2\right)\frac{\partial^2 \xi}{\partial x^2}
+a^2_{\infty}\frac{\partial^2 \xi}{\partial z^2}
-\lambda^2 \xi-2M\lambda a_{\infty}\frac{\partial \xi}{\partial x}=0.
\end{equation}
\indent
To describe the general solution of equation (\ref{LaplaceEquation}) satisfying
the far field condition we consider, following \cite{Ba2}
$$D(\omega,\lambda)=M^2\left(\frac{\lambda}{U}\right)^2+2i\frac{\lambda}{U}M^2\omega
+\left(1-M^2\right)\omega^2,$$
and prove two lemmas below.\\

\begin{lemma}\label{RootD}\ There exists a function $\sqrt{D(\omega,\lambda)}$, analytic
with respect to complex variable $\frac{\lambda}{U}+i\omega\ (\frac{\lambda}{U}
\in\C,\ \omega\in \R)$ in the half-plane ${\dis \mbox{Re}\lambda>\sigma_a}$,
and such that $\mbox{Re}\sqrt{D(\omega,\lambda)}>0$.
\end{lemma}
\indent
{\bf Proof.}\ Representing $D(\omega,\lambda)$ as
$$D(\omega,\lambda)=M^2\left(\frac{\lambda}{U}\right)^2+2i\frac{\lambda}{U}M^2\omega
+\left(1-M^2\right)\omega^2=M^2\left(\frac{\lambda}{U}+i\omega\right)^2+\omega^2,$$
we obtain that the image of the half-plane
${\dis \mbox{Re}\lambda>\sigma_a}$
under the map $D(\omega,\lambda)$ is contained in the domain
$\C\setminus{\R^{-}}$. Then the branch of the function $\sqrt{}$ considered
on the complex plane with the cut along the negative part of the real axis is well
defined and analytic on the image of $D$, and its real part satisfies condition
of the Lemma. Therefore, the composition $\sqrt{D}$ is also analytic and
satisfies the same condition.
\qed

\begin{lemma}\label{Fourier}\ The following equality holds
$$\frac{e^{d(\lambda)x}}{\sqrt{1-M^2}}
K_0\left(r(\lambda)\left(\frac{x^2}{1-M^2}+z^2\right)^{\frac{1}{2}}\right)
={\cal F}\left[\frac{e^{\dis -z\sqrt{D(\omega,\lambda)} }}
{2\sqrt{D(\omega,\lambda)}}\right],$$
where ${\cal F}$ denotes the Fourier transform, or
$$\frac{e^{d(\lambda)x}}{\sqrt{1-M^2}}
K_0\left(r(\lambda)\left(\frac{x^2}{1-M^2}+z^2\right)^{\frac{1}{2}}\right)$$
$$=\int_{-\infty}^\infty e^{ix\omega}
\frac{e^{\dis -z\left((1-M^2)(\omega+id(\lambda))^2+r^2(\lambda)\right)^{\frac{1}{2}}} }
{2\sqrt{(1-M^2)(\omega+id(\lambda))^2+r^2(\lambda)}}d\omega.$$
\end{lemma}
\indent
{\bf Proof.}\ First, we represent $D\left(\omega,\lambda\right)$ as
$$D\left(\omega,\lambda\right)=\left(1-M^2\right)\omega^2+2i\frac{\lambda}{U}M^2\omega
+M^2\left(\frac{\lambda}{U}\right)^2$$
$$=\left(\omega\sqrt{1-M^2}+i\frac{\lambda M^2}{U\sqrt{1-M^2}}
\right)^2+\left(\frac{\lambda M^2}{U\sqrt{1-M^2}}\right)^2
+M^2\left(\frac{\lambda}{U}\right)^2$$
$$=(1-M^2)(\omega+id(\lambda))^2+r^2(\lambda).$$
\indent
Changing variables in equality (\cite{EMOT})
$$K_0\left(r(x^2+z^2)^{\frac{1}{2}}\right)
=\int_{-\infty}^\infty e^{ixu}
\frac{e^{\dis -z\left(u^2+r^2\right)^{\frac{1}{2}}} }{2\sqrt{u^2+r^2}}du,$$
we obtain
$$K_0\left(r(x^2+z^2)^{\frac{1}{2}}\right)
=\int_{-\infty}^\infty e^{ix\sqrt{1-M^2}\omega}
\frac{e^{\dis -z\left((1-M^2)\omega^2+r^2\right)^{\frac{1}{2}}} }
{2\sqrt{(1-M^2)\omega^2+r^2}}d\left(\sqrt{1-M^2}\omega\right),$$
and then
$$\frac{1}{\sqrt{1-M^2}}K_0\left(r\left(\frac{x^2}{1-M^2}+z^2\right)^{\frac{1}{2}}\right)
=\int_{-\infty}^\infty e^{ix\omega}
\frac{e^{\dis -z\left((1-M^2)\omega^2+r^2\right)^{\frac{1}{2}}} }
{2\sqrt{(1-M^2)\omega^2+r^2}}d\omega.$$
\indent
We transform the equality above by integrating function
$$g(x,w)=e^{ixw}\frac{e^{\dis -z\left((1-M^2)w^2+r^2\right)^{\frac{1}{2}}} }
{\sqrt{(1-M^2)w^2+r^2}},\hspace{0.1in}w\in\C,$$
analytic with respect to $w$, over the piecewise linear contour
$$\left[-C, C, C+id, -C+id\right]\in \C,\hspace{0.1in}\mbox{with}\hspace{0.05in}
C\in\R, C>0,\ d\in\C,\ \mbox{Re}d>0.$$
Then we obtain
\begin{equation}\label{ContourIntegral}
\int_{-C}^Cg(x,w)dw+\int_{C}^{C+id}g(x,w)dw+\int_{C+id}^{-C+id}g(x,w)dw
+\int_{-C+id}^{-C}g(x,w)dw=0.
\end{equation}
\indent
For $C$ large enough we have the following estimates for $w=u+iv \in \left[-C, -C+id\right]$,
and $w \in \left[C, C+id\right]$
$$\left|e^{ix(u+iv)}\right|<e^{|x|\cdot\mbox{\scriptsize Re}d},
\hspace{0.1in}\left|\sqrt{(1-M^2)w^2+r^2}\right|>\sqrt{1-M^2}\frac{C}{2},$$
$$\left|e^{\dis -z\left((1-M^2)w^2+r^2\right)^{\frac{1}{2}}}\right|
<e^{-z\sqrt{1-M^2}\frac{C}{2}},$$
and therefore for $z>0$
$$\left|\int_{C}^{C+id}g(x,w)dw\right|,\left|\int_{-C}^{-C+id}g(x,w)dw\right|\to 0
\hspace{0.1in}\mbox{as}\hspace{0.05in}C\to\infty.$$
\indent
Using the last estimate in (\ref{ContourIntegral}) we obtain equality
$$\int_{-\infty}^\infty e^{ix\omega}
\frac{e^{\dis -z\left((1-M^2)\omega^2+r^2(\lambda)\right)^{\frac{1}{2}}} }
{\sqrt{(1-M^2)\omega^2+r^2(\lambda)}}d\omega$$
$$=\int_{-\infty}^\infty e^{ix(\omega+id(\lambda))}
\frac{e^{\dis -z\left((1-M^2)(\omega+id(\lambda))^2+r^2(\lambda)\right)^{\frac{1}{2}}} }
{\sqrt{(1-M^2)(\omega+id(\lambda))^2+r^2(\lambda)}}d(\omega+id(\lambda)),$$
and finally
\begin{equation}\label{SFormula}
\frac{e^{d(\lambda)x}}{\sqrt{1-M^2}}
K_0\left(r(\lambda)\left(\frac{x^2}{1-M^2}+z^2\right)^{\frac{1}{2}}\right)
\end{equation}
$$=\int_{-\infty}^\infty e^{ix\omega}
\frac{e^{\dis -z\left((1-M^2)(\omega+id(\lambda))^2+r^2(\lambda)\right)^{\frac{1}{2}}} }
{2\sqrt{(1-M^2)(\omega+id(\lambda))^2+r^2(\lambda)}}d\omega.$$
\qed

\indent
Using now Lemmas ~\ref{RootD} and ~\ref{Fourier} we consider a special representation
of the general solution of (\ref{LaplaceEquation}). Namely, using notations of
Lemma~\ref{Fourier}, and denoting
$$S(x,z,\lambda)=-\frac{e^{d(\lambda)x}}{\sqrt{1-M^2}}
K_0\left(r(\lambda)\left(\frac{x^2}{1-M^2}+z^2\right)^{\frac{1}{2}}\right),$$
we consider
\begin{equation}\label{LapSolution}
\begin{array}{ll}
{\dis \xi(x,z,\lambda)=\int_{-\infty}^\infty S(x-y,z,\lambda)v_a(y,\lambda)dy }
\vspace{0.1in}\\
{\dis =-\frac{e^{d(\lambda)x}}{\sqrt{1-M^2}}\int_{-\infty}^\infty e^{-d(\lambda)y}
K_0\left(r(\lambda)\left(\frac{(x-y)^2}{1-M^2}+z^2\right)^{\frac{1}{2}}\right)v_a(y,\lambda)dy. }
\end{array}
\end{equation}

\begin{proposition}\label{LaplaceSolution}\ Function $\xi$ defined by formula
(\ref{LapSolution}) satisfies equation (\ref{LaplaceEquation}).
If
$$\xi(x,z,\lambda),\ \frac{\partial^2 \xi}{\partial x^2}(x,z,\lambda),\
\frac{\partial^2 \xi}{\partial z^2}(x,z,\lambda),\ |\eta|^2\xi(x,z,\lambda),\
|\eta|\frac{\partial\xi(x,z,\lambda)}{\partial x}\in L^1(\R_{\eta}),$$
where $\lambda=\sigma+i\eta$,
then the inverse Laplace transform of $\xi$, defined by the formula (\cite{Boc})
\begin{equation}\label{InverseLaplace}
\phi(x,z,t)=\frac{1}{2\pi}\int_{\sigma-i\infty}^{\sigma+i\infty}e^{(\sigma+i\eta)t}
\xi(x,z,\sigma+i\eta)d\eta
\end{equation}
satisfies equation (\ref{LinearEquation}).
\end{proposition}
\indent
{\bf Proof.}\  To prove that $\xi$ defined above satisfies
equation (\ref{LaplaceEquation})
it suffices to prove that function $S$ satisfies the same equation. For $S$ we have
$$a^2_{\infty}\left(1-M^2\right)\frac{\partial^2 S}{\partial x^2}
+a^2_{\infty}\frac{\partial^2 S}{\partial z^2}
-\lambda^2 S-2M\lambda a_{\infty}\frac{\partial S}{\partial x}$$
$$=a^2_{\infty}\left[\frac{\partial^2 S}{\partial z^2}
+\left(1-M^2\right)\frac{\partial^2 S}{\partial x^2}
-M^2\left(\frac{\lambda}{U}\right)^2 S-2M^2\frac{\lambda}{U}\frac{\partial S}{\partial x}\right].$$
\indent
Using then formula (\ref{SFormula}), we obtain
$$\frac{\partial^2 S}{\partial z^2}
+\left(1-M^2\right)\frac{\partial^2 S}{\partial x^2}
-M^2\left(\frac{\lambda}{U}\right)^2 S-2M^2\frac{\lambda}{U}\frac{\partial S}{\partial x}$$
$$=-\int_{-\infty}^\infty e^{ix\omega}\left((1-M^2)(\omega+id)^2+r^2\right)
\frac{e^{\dis -z\left((1-M^2)(\omega+id)^2+r^2\right)^{\frac{1}{2}}} }
{2\sqrt{(1-M^2)(\omega+id)^2+r^2}}d\omega$$
$$+\int_{-\infty}^\infty e^{ix\omega}\left(1-M^2\right)\omega^2
\frac{e^{\dis -z\left((1-M^2)(\omega+id)^2+r^2\right)^{\frac{1}{2}}} }
{2\sqrt{(1-M^2)(\omega+id)^2+r^2}}d\omega$$
$$+\int_{-\infty}^\infty e^{ix\omega}M^2\left(\frac{\lambda}{U}\right)^2
\frac{e^{\dis -z\left((1-M^2)(\omega+id)^2+r^2\right)^{\frac{1}{2}}} }
{2\sqrt{(1-M^2)(\omega+id)^2+r^2}}d\omega$$
$$+\int_{-\infty}^\infty e^{ix\omega}2M^2\frac{\lambda}{U}i\omega
\frac{e^{\dis -z\left((1-M^2)(\omega+id)^2+r^2\right)^{\frac{1}{2}}} }
{2\sqrt{(1-M^2)(\omega+id)^2+r^2}}d\omega=0.$$
\indent
To prove that function $\phi$ defined by formula (\ref{InverseLaplace}) satisfies
equation (\ref{LinearEquation}) we apply the inverse Laplace transform to equality
$$a^2_{\infty}\left(1-M^2\right)\frac{\partial^2 \xi}{\partial x^2}
+a^2_{\infty}\frac{\partial^2 \xi}{\partial z^2}
-\lambda^2 \xi-2M\lambda a_{\infty}\frac{\partial
\xi}{\partial x}=0$$
and obtain equation (\ref{LinearEquation}) for $\phi$.
\qed

\section{Boundary Conditions.}
\label{Boundary}

\indent
In this section we reformulate the boundary conditions of section ~\ref{Introduction}
in terms of function $v_a(y,\lambda)$ from formula (\ref{LapSolution}).\\
\indent
To check the flow tangency condition (\ref{Flow_Tangency}) we use formulas (\ref{LapSolution})
and (\ref{SFormula}), and obtain
\begin{equation}\label{inversion}
\frac{\partial}{\partial z}\xi(x,z,\lambda)\Big|_{z=0}
=\frac{\partial}{\partial z}\int_{-\infty}^\infty
S(x-y,z,\lambda)v_a(y,\lambda)dy\Big|_{z=0}
\end{equation}
$$=-\frac{\partial}{\partial z}\int_{-\infty}^\infty v_a(y,\lambda)dy
\int_{-\infty}^\infty e^{i(x-y)\omega}
\frac{e^{\dis -z\left((1-M^2)(\omega+id(\lambda))^2+r^2(\lambda)\right)^{\frac{1}{2}}} }
{2\sqrt{(1-M^2)(\omega+id(\lambda))^2+r^2(\lambda)}}d\omega\Big|_{z=0}$$
$$=\frac{1}{2}\int_{-\infty}^\infty e^{ix\omega}d\omega\int_{-\infty}^\infty
e^{-iy\omega}v_a(y,\lambda)dy=\pi\cdot v_a(x,\lambda),$$
which, after comparison with equality (\ref{Flow_Tangency}) leads to a unique choice
\begin{equation}\label{v=w}
v_a(x,\lambda)=\frac{1}{\pi}\widehat{w_a}(x,\lambda)\ \mbox{for}\ |x|<1.
\end{equation}
\indent
To satisfy the Kutta-Joukowski boundary condition (\ref{Kutta}) we should have
$${\dis \left(\frac{\partial\phi}{\partial t}
+U\frac{\partial\phi}{\partial x}\right)\Big|_{z=0}
=0\hspace{0.1in}\mbox{for}\ 1<|x|<A,}$$
or equality
$$\lambda\xi(x,0,\lambda)
+U\frac{\partial\xi}{\partial x}(x,0,\lambda)
=0\hspace{0.1in}\mbox{for}\ 1<|x|<A$$
for function $\xi$.\\
\indent
Substituting $\xi$ from formula (\ref{LapSolution}) into equality above
we obtain the following condition for $1<|x|<A$:
\begin{equation}\label{ZeroCondition}
\begin{array}{ll}
{\dis 0=\left(\lambda+U\frac{\partial}{\partial x}\right)\xi(x,0,\lambda) }
\vspace{0.1in}\\
{\dis =-\left(\lambda+U\frac{\partial}{\partial x}\right)
\frac{e^{d(\lambda)x}}{\sqrt{1-M^2}}\int_{-\infty}^\infty e^{-d(\lambda)y}
K_0\left(r(\lambda)\left(\frac{(x-y)^2}{1-M^2}+z^2\right)^{\frac{1}{2}}\right)
v_a(y,\lambda)dy\Big|_{z=0}. }
\end{array}
\end{equation}
\indent
To reformulate the last condition as an integral equation we
use condition (\ref{v=w}), and define
$$g_a(x,\lambda)=\frac{e^{d(\lambda)x}}{\pi}\int_{-1}^1 e^{-d(\lambda)y}
R(x-y,\lambda)\widehat{w_a}(y,\lambda)dy\ \mbox{for}\ 1<|x|<A,$$
with kernel $R(x,\lambda)$ defined by the formula
\begin{equation}\label{RKernel}
R(x,\lambda)=\left[\left(\lambda+Ud(\lambda)\right)K_0
\left(\frac{r(\lambda)|x|}{\sqrt{1-M^2}}\right)
+U\frac{\partial}{\partial x}K_0\left(\frac{r(\lambda)|x|}{\sqrt{1-M^2}}\right)\right].
\end{equation}
Then condition (\ref{ZeroCondition}) will be satisfied if $v_a$ will
satisfy the following integral equation
$$e^{d(\lambda)x}\int_{|y|>1} e^{-d(\lambda)y}
R(x-y,\lambda)v_a(y,\lambda)dy=-g_a(x,\lambda)\ \mbox{for}\ 1<|x|<A.$$
\indent
Further simplifying the equation above we consider
$h_a(y,\lambda):=e^{-d(\lambda)y}\cdot v_a(y,\lambda)$ as an
unknown function, and rewrite it as
\begin{equation}\label{MainEquation}
\int_{|y|>1}R(x-y,\lambda)h_a(y,\lambda)dy
=f_a(x,\lambda)\hspace{0.1in}\mbox{for}\ 1<|x|<A,
\end{equation}
where $f_a(x,\lambda)=-e^{-d(\lambda)x}\cdot\chi_A(x) g_a(x,\lambda)$
is defined for
$$\left\{(x,\lambda)\in \R\times\C:
|x|>1, \mbox{Re}\lambda\in[\sigma_1,\sigma_2]\right\}$$
by the formula
\begin{equation}\label{fDefinition}
f_a(x,\lambda)=-\frac{\chi_A(x)}{\pi}\int_{-1}^1 e^{-d(\lambda)y}
R(x-y,\lambda)\widehat{w_a}(y,\lambda)dy
\end{equation}
with
$$\chi_A(x)=\left\{\begin{array}{ll}
1\hspace{0.1in}\mbox{if}\ x\in [-A,A]\setminus[-1,1],\vspace{0.1in}\\
0\hspace{0.1in}\mbox{otherwise}.
\end{array}\right.$$

\section{Cofinite Hilbert transform.}
\label{CofiniteTransform}

\indent
As a first step in the analysis of solvability of (\ref{MainEquation}) we prove
solvability for the operator, closely related to operator
${\cal R}_{\lambda}$ from (\ref{MainEquation}), and which in analogy with the Tricomi's
definition of the finite Hilbert transform \cite{Tr} we call the
{\it cofinite Hilbert transform}.\\
\indent
We define the cofinite Hilbert transform on the set of functions on
$$I^c(1)=\R\setminus[-1,1]$$
by the formula
\begin{equation}\label{POperator}
{\cal P}[h](x)=\frac{1}{\pi}\int_{|y|>1}\frac{h(y)}{y-x}dy\
\mbox{for}\ |x|>1,
\end{equation}
where the integral
$$\int_{|y|>1}=\int_{-\infty}^{-1}+\int_1^{\infty}$$
is understood in the sense of Cauchy's principal value.\\
\indent
In the proposition below we prove solvability for the nonhomogeneous integral
equation with operator ${\cal P}$ in weighted spaces
$${\cal L}^p\left(I^c(1)\right)=\left\{f:\int_{|x|>1}|x|^{p-2}\left|f(x)\right|^pdx
<\infty\right\}$$
with
$$\|f\|_{{\cal L}^p\left(I^c(1)\right)}
=\left(\int_{|x|>1}|x|^{p-2}\left|f(x)\right|^pdx\right)^{1/p}.$$

\begin{proposition}\label{Cofinite}\ For any function $f \in {\cal L}^q\left(I^c(1)\right)$
with ${\dis q>\frac{4}{3} }$ there exists a solution $h$ of equation
\begin{equation}\label{Pequation}
{\cal P}[h]=f,
\end{equation}
such that $h \in {\cal L}^p\left(I^c(1)\right)$ for any ${\dis p<\frac{4}{3} }$.
\end{proposition}
\indent
{\bf Proof.}\ We consider the following diagram of transformations
\begin{equation}\label{Diagram}
\begin{array}{ccc}
L^p\left(I(1)\right)&\mapright{-\cal T}&L^p\left(I(1)\right)
\vspace{0.1in}\\
\mapdown{\Theta}&&\mapdown{\Theta}\vspace{0.1in}\\
{\cal L}^p\left(I^c(1)\right)&\mapright{\cal P}&{\cal L}^p\left(I^c(1)\right),
\end{array}
\end{equation}
where ${\cal T}$ is the finite Hilbert transform, ${\cal P}$ is the
cofinite Hilbert transform, and
$$\Theta:L^p\left(I(1)\right)\to {\cal L}^p\left(I^c(1)\right)$$
is defined by the formula
\begin{equation}\label{ThetaFormula}
\Theta[f](x)=\frac{1}{x}f\left(\frac{1}{x}\right).
\end{equation}
\indent
To prove that the maps in diagram (\ref{Diagram}) are well defined we use
equality
$$\left\|\Theta[f]\right\|_{{\cal L}^p}^p
=\int_{|x|>1}|x|^{p-2}\left|\Theta[f](x)\right|^pdx
=\int_{|x|>1}|x|^{p-2}\frac{\left|f\left(\frac{1}{x}\right)\right|^p}{|x|^p}dx$$
$$=-\int_1^{-1}\left|f\left(t\right)\right|^pdt
=\|f\|_p^p,$$
and notice that for
$$\Theta^*:{\cal L}^p\left(I^c(1)\right)\to L^p\left(I(1)\right)$$
defined by the same formula
$$\Theta^*[f](x)=\frac{1}{x}f\left(\frac{1}{x}\right),$$
we have
\begin{equation}\label{Thetaand*}
\Theta\circ\Theta^*\left[f\right](x)=\Theta
\left[\frac{1}{y}f\left(\frac{1}{y}\right)\right](x)
=\frac{1}{x}\cdot xf(x)=f(x).
\end{equation}
\indent
Diagram (\ref{Diagram}) is commutative, as can be seen from equality
$${\cal P}\left[\Theta[f]\right](x)=\frac{1}{\pi}\int_{|y|>1}
\frac{f\left(\frac{1}{y}\right)}{y(y-x)}dy
=\frac{1}{\pi}\int_{-1}^1t\frac{f\left(t\right)}{(\frac{1}{t}-x)t^2}dt$$
$$=\frac{1}{\pi x}\int_{-1}^1\frac{f\left(t\right)}{\frac{1}{x}-t}dt
=\Theta\left[-{\cal T}[f]\right].$$
\indent
To "invert" operator ${\cal P}$ we use commutativity of diagram (\ref{Diagram}),
relation (\ref{Thetaand*}), and operator (\cite{So},\cite{Tr})
$${\cal T}^{-1}:L^{\frac{4}{3}+}\left(I(1)\right)\to L^{\frac{4}{3}-}\left(I(1)\right),$$
defined by the formula
$${\cal T}^{-1}[g](x)=-\frac{1}{\pi}\int_{-1}^{1}\sqrt{\frac{1-y^2}{1-x^2}}
\frac{g(y)}{y-x}dy,$$
and satisfying
$${\cal T}\circ{\cal T}^{-1}[f]=f.$$
\indent
Namely, we define operator
$${\cal P}^{-1}:{\cal L}^{\frac{4}{3}+}\left(I^c(1)\right)\to
{\cal L}^{\frac{4}{3}-}\left(I^c(1)\right)$$
by the formula
$${\cal P}^{-1}[f]=-\Theta\circ{\cal T}^{-1}\circ\Theta^*[f].$$
\indent
Then
$${\cal P}\circ{\cal P}^{-1}[f]=-{\cal P}\circ\Theta\circ{\cal T}^{-1}\circ\Theta^*[f]
=\Theta\circ{\cal T}\circ{\cal T}^{-1}\circ\Theta^*[f]
=\Theta\circ\Theta^*[f]=f,$$
and we obtain the statement of the proposition for
$$h={\cal P}^{-1}[f].$$
\indent
To find an explicit formula for ${\cal P}^{-1}$ we use explicit formulas for $\Theta$
and ${\cal T}^{-1}$, and obtain
\begin{equation}\label{P-1Forfula}
\begin{array}{ll}
{\dis {\cal P}^{-1}[f](x)=\frac{1}{\pi x}\int_{-1}^{1}\sqrt{\frac{1-y^2}{1-1/x^2}}
\cdot\frac{f(1/y)}{y(y-1/x)}dy }\vspace{0.1in}\\
{\dis =\frac{|x|}{\pi}\int_{-1}^{1}\sqrt{\frac{1-y^2}{x^2-1}}
\left[\frac{1}{y}f\left(\frac{1}{y}\right)\right]\frac{dy}{xy-1}. }
\end{array}
\end{equation}
\qed

\indent
{\bf Remark.}\hspace{0.1in}Following \cite{Tr} we notice that solution of equation
(\ref{Pequation}) is unique in ${\cal L}^2\left(I^c(1)\right)$, but is not unique in larger
spaces. Namely, function
$$h(x)=\frac{1}{\sqrt{x^2-1}}$$
is a solution, and the only one in ${\cal L}^{2-}\left(I^c(1)\right)$ up to the
linear dependence, of the homogeneous equation
$${\cal P}[h]=0.$$
\qed

\section{Solvability of equation (\ref{MainEquation}).}
\label{Solvability}

\indent
From the asymptotic expansions of $K_0(\zeta)$
(see \cite{EMOT}, \cite{GR}) we obtain the following representations of the function
$R(x,\lambda)$ for $\lambda$ such that $\mbox{Re}\lambda\in[\sigma_1,\sigma_2]$
with $\sigma_1>\sigma_a$:
\begin{equation}\label{RRepresentation}
\begin{array}{ll}
{\dis R(x,\lambda)=-\frac{U}{x}+\lambda\log\left(\lambda|x|\right)\alpha(\lambda|x|)
+\lambda\beta\left(\lambda|x|\right)+\gamma\left(\lambda|x|\right)
\hspace{0.1in}\mbox{for}\ |\lambda x|\leq B, }\vspace{0.1in}\\
{\dis R(x,\lambda)=\lambda\delta\left(\lambda|x|\right)\frac{e^{-(\sigma+i\eta)|x|}}
{\sqrt{|\lambda|\cdot|x|}}\hspace{0.1in}\mbox{for}\ |\lambda x|> B, }
\end{array}
\end{equation}
where $\alpha(\zeta)$, $\beta(\zeta)$, $\gamma(\zeta)$, and $\delta(\zeta)$ are bounded
analytic functions on $\mbox{Re}\zeta>\epsilon>0$ and $B>0$ is some constant.\\
\indent
Using representations (\ref{RRepresentation}) we introduce function $M(x,\lambda)$,
analytic with respect to $\lambda\in \left\{\mbox{Re}\lambda>\sigma_a\right\}$,
uniquely defined by (\ref{RRepresentation}), and such that
$${\dis R(x,\lambda)=-\frac{U}{x}+M(x,\lambda). }$$
We consider then operators
$${\cal M}_{\lambda}[f](x)=\int_{|y|>1}\chi_A(x)M(x-y,\lambda)f(y)dy,$$
and
$${\cal R}_{\lambda}=\pi U\cdot{\cal P}+{\cal M}_{\lambda}.$$
\indent
In the next proposition we prove compactness of the operator
${\dis \frac{1}{\pi U}{\cal M}_{\lambda}\circ{\cal P}^{-1} }$
on ${\cal L}^{2}\left(I^c(1)\right)$.\\

\begin{proposition}\label{RFredholm}\ For any fixed $\lambda\in\C$ operator
${\dis {\cal N}_{\lambda}=\frac{1}{\pi U}{\cal M}_{\lambda}\circ{\cal P}^{-1} }$
is compact on ${\cal L}^{2}\left(I^c(1)\right)$, and therefore operator
\begin{equation}\label{GOperator}
{\cal G}_{\lambda}={\cal R}_{\lambda}\circ{\cal P}^{-1}
=\left(\pi U\cdot{\cal P}+{\cal M}_{\lambda}\right)\circ{\cal P}^{-1}
=\pi U\left({\cal I}+{\cal N}_{\lambda}\right),
\end{equation}
where ${\cal I}$ is the identity operator, is a Fredholm operator on
${\cal L}^{2}\left(I^c(1)\right)=L^{2}\left(I^c(1)\right)$. In addition, kernel
$N(x,y,\lambda)$ of the operator ${\cal N}_{\lambda}$ admits estimate
\begin{equation}\label{NEstimate}
\int_{\R^2}|N(x,y,\lambda)|^2dxdy<C|\lambda\log{\lambda}|^2
\end{equation}
with constant $C$ independent of $\lambda$.
\end{proposition}
\indent
{\bf Proof.}\ Using formula (\ref{P-1Forfula}) for ${\cal P}^{-1}$, we obtain
$${\cal N}_{\lambda}[g](x)={\cal M}_{\lambda}
\left[\frac{|x|}{\pi^2 U}\int_{-1}^{1}\sqrt{\frac{1-u^2}{x^2-1}}
\left[\frac{1}{u}g\left(\frac{1}{u}\right)\right]\frac{du}{xu-1}\right]$$
$$={\cal M}_{\lambda}\left[\frac{|x|}{\pi^2 U}\int_{|y|>1}\sqrt{\frac{1-\frac{1}{y^2}}{x^2-1}}
y^2g(y)\frac{dy}{y^2(x-y)}\right]$$
$$=\frac{\chi_A(x)}{\pi^2 U}\int_{|u|>1}M(x-u,\lambda)du
\int_{|y|>1}\frac{|u|\sqrt{y^2-1}}{|y|\sqrt{u^2-1}}g(y)\frac{dy}{(u-y)}$$
$$=\frac{\chi_A(x)}{\pi^2 U}\int_{|y|>1}g(y)dy
\int_{|u|>1}M(x-u,\lambda)\frac{|u|\sqrt{y^2-1}}{|y|\sqrt{u^2-1}}
\frac{du}{(u-y)}=\int_{|y|>1}N(x,y,\lambda)g(y)dy$$
with kernel
$$N(x,y,\lambda)=\frac{\chi_A(x)}{\pi^2 U}
\int_{|u|>1}M(x-u,\lambda)\frac{|u|\sqrt{y^2-1}}{|y|\sqrt{u^2-1}}
\frac{du}{(u-y)}.$$
\indent
To prove compactness of operator ${\cal N}_{\lambda}$ we use representation
$$N(x,y,\lambda)=\frac{1}{\pi^2 U}\left[N_1(x,y,\lambda)+N_2(x,y,\lambda)\right],$$
with
$$N_1(x,y,\lambda)=\frac{\chi_A(x)}{|y|}\int_{|u|>1}M(x-u,\lambda)|u|
\frac{du}{(u-y)},$$
and
$$N_2(x,y,\lambda)=\frac{\chi_A(x)}{|y|}\int_{|u|>1}M(x-u,\lambda)
\frac{|u|\left(\sqrt{y^2-1}-\sqrt{u^2-1}\right)}{\sqrt{u^2-1}}
\frac{du}{(u-y)}$$
$$=-\frac{\chi_A(x)}{|y|}\int_{|u|>1}M(x-u,\lambda)
\frac{|u|\left(y+u\right)du}{\left(\sqrt{y^2-1}+\sqrt{u^2-1}\right)\sqrt{u^2-1}},$$
and prove Hilbert-Schmidt property (cf.\cite{L})
of kernels $N_1(x,y,\lambda)$ and $N_2(x,y,\lambda)$.\\
\indent
For $N_1(x,y,\lambda)$ we notice that for fixed $x$ satisfying $1<|x|<A$
$$\int_{|u|>1}M(x-u,\lambda)|u|\frac{du}{(u-y)}$$
is a multiple of the Hilbert transform of an $L^{2}\left(I^c(1)\right)$ - function
$M(x-u,\lambda)|u|$ with
$$\|M(x-u,\lambda)|u|\|_{L^{2}_u}<\infty.$$
Therefore we have
\begin{equation}\label{N1Estimate}
\begin{array}{lll}
{\dis \int_{1<|x|<A}dx\int_{|y|>1}dy\left|N_1(x,y,\lambda)\right|^2 }\vspace{0.1in}\\
{\dis =\int_{1<|x|<A}dx\int_{|y|>1}dy\left|\frac{1}{|y|}\int_{|u|>1}M(x-u,\lambda)|u|
\frac{du}{(u-y)}\right|^2 }\vspace{0.1in}\\
{\dis < C\int_{1<|x|<A}dx\left\|M(x-u,\lambda)|u|\right\|_{L^{2}_u}^2<\infty. }
\end{array}
\end{equation}
\indent
For $N_2(x,y,\lambda)$ we have
\begin{equation}\label{N2Integral}
\begin{array}{llll}
{\dis \int_{1<|x|<A}dx\int_{|y|>1}dy\left|N_2(x,y,\lambda)\right|^2 }\vspace{0.1in}\\
{\dis =\int_{1<|x|<A}dx\int_{|y|>1}\frac{dy}{|y|^2}\left|\int_{|u|>1}M(x-u,\lambda)
\frac{|u|\left(y+u\right)du}{\left(\sqrt{y^2-1}+\sqrt{u^2-1}\right)\sqrt{u^2-1}}\right|^2 }
\vspace{0.1in}\\
{\dis \leq 2\int_{1<|x|<A}dx\int_{|y|>1}\frac{dy}{|y|^2}
\left|\int_0^{\infty}M(x-\sqrt{t^2+1},\lambda)
\frac{\left(y+\sqrt{t^2+1}\right)dt}{\left(\sqrt{y^2-1}+t\right)}\right|^2 }
\vspace{0.1in}\\
{\dis +2\int_{1<|x|<A}dx\int_{|y|>1}\frac{dy}{|y|^2}
\left|\int_0^{\infty}M(x+\sqrt{t^2+1},\lambda)
\frac{\left(y-\sqrt{t^2+1}\right)dt}{\left(\sqrt{y^2-1}+t\right)}\right|^2, }
\end{array}
\end{equation}
where we changed variable to $t=\sqrt{u^2-1}$.\\
\indent
Both integrals of the right hand side of (\ref{N2Integral}) are estimated analogously,
therefore we will present an estimate of the first of them only.\\
\indent
For $1<|x|<A$ and $|y|>2$ we have inequality
\begin{equation}\label{y-tEstimate}
\left|\frac{y+\sqrt{t^2+1}}{\sqrt{y^2-1}+t}\right|<C
\end{equation}
for some $C$ independent of $y$, and therefore, using representations (\ref{RRepresentation}),
we obtain
\begin{equation}\label{N2Estimate1}
\begin{array}{ll}
{\dis \int_{1<|x|<A}dx\int_{|y|>2}\frac{dy}{|y|^2}
\left|\int_0^{\infty}M(x-\sqrt{t^2+1},\lambda)
\frac{\left(y+\sqrt{t^2+1}\right)dt}{\left(\sqrt{y^2-1}+t\right)}\right|^2 }\vspace{0.1in}\\
{\dis \leq C^2\int_{1<|x|<A}dx\int_{|y|>2}\frac{dy}{|y|^2}
\left|\int_0^{\infty}M(x-\sqrt{t^2+1},\lambda)dt\right|^2<\infty. }
\end{array}
\end{equation}
\indent
For $1<|x|<A$, $1<|y|<2$, and $t>A+B$ we again use inequality (\ref{y-tEstimate}) and obtain
\begin{equation}\label{N2Estimate2}
\begin{array}{ll}
{\dis \int_{1<|x|<A}dx\int_{1<|y|<2}\frac{dy}{|y|^2}
\left|\int_{A+B}^{\infty}M(x-\sqrt{t^2+1},\lambda)
\frac{\left(y+\sqrt{t^2+1}\right)dt}{\left(\sqrt{y^2-1}+t\right)}\right|^2 }\vspace{0.1in}\\
{\dis \leq C^2\int_{1<|x|<A}dx\int_{|y|>2}\frac{dy}{|y|^2}
\left|\int_0^{\infty}M(x-\sqrt{t^2+1},\lambda)dt\right|^2<\infty. }
\end{array}
\end{equation}
\indent
For $1<|x|<A$, $1<|y|<2$, and $t<A+B$ we have
$$\int_{1<|x|<A}dx\int_{1<|y|<2}\frac{dy}{|y|^2}\left|\int_0^{A+B}M(x-\sqrt{t^2+1},\lambda)
\frac{\left(y+\sqrt{t^2+1}\right)dt}{\left(\sqrt{y^2-1}+t\right)}\right|^2$$
$$< C\int_{1<|x|<A}dx\int_{1<|y|<2}dy\left|\int_0^{A+B}
\frac{M(x-\sqrt{t^2+1},\lambda)dt}{\left(\sqrt{y^2-1}+t\right)}\right|^2$$
$$< C|\lambda\log{\lambda}|^2\int_{1<|x|<A}dx\int_{1<|y|<2}dy\left|\int_0^{A+B}
\frac{dt}{\left(\sqrt{y^2-1}+t\right)}\right|^2$$
$$+ C|\lambda|^2\int_{1<|x|<A}dx\int_{1<|y|<2}dy\left|\int_0^{A+B}
\frac{\log{|x-\sqrt{t^2+1}|}}{\left(\sqrt{y^2-1}+t\right)}dt\right|^2$$
$$<C|\lambda|^2\left(|\log{\lambda}|^2+\int_{1<|x|<A}dx\int_{1<|y|<2}dy\left|\int_0^{A+B}
\frac{\log{|x-\sqrt{t^2+1}|}}{\left(\sqrt{y^2-1}+t\right)}dt\right|^2\right),$$
where we used representation
$$M(x-\sqrt{t^2+1},\lambda)=\lambda\left(\log{\lambda}+\log{|x-\sqrt{t^2+1}|}\right)
\alpha(\lambda|x-\sqrt{t^2+1}|)$$
$$+\lambda \beta(\lambda|x-\sqrt{t^2+1}|)+\gamma(\lambda|x-\sqrt{t^2+1}|)$$
for $1<|x|<A$ and $0<t<A+B$, which is a corollary of (\ref{RRepresentation}).\\
\indent
To estimate the last integral we represent it as
$$\int_{1<|x|<A}dx\int_{1<|y|<2}dy\left|\int_0^{A+B}
\frac{\log{|x-\sqrt{t^2+1}|}}{\left(\sqrt{y^2-1}+t\right)}dt\right|^2$$
$$=\int_{1<|x|<A}dx\int_{1<|y|<2}dy\left|\int_{S(x,y)}
\frac{\log{|x-\sqrt{t^2+1}|}}{\left(\sqrt{y^2-1}+t\right)}dt\right|^2$$
$$+\int_{1<|x|<A}dx\int_{1<|y|<2}dy
\left|\int_{\left\{\left[0,A+B\right]\setminus S(x,y)\right\}}
\frac{\log{|x-\sqrt{t^2+1}|}}{\left(\sqrt{y^2-1}+t\right)}dt\right|^2$$
where $S(x,y)=\left\{t:\left|x-\sqrt{t^2+1}\right|
\geq\frac{1}{2}\left(x-1\right)\sqrt{y^2-1}\right\}$.\\
\indent
Then for $S(x,y)$ we have
$$\int_{1<|x|<A}dx\int_{1<|y|<2}dy\left|\int_{S(x,y)}
\frac{\log{|x-\sqrt{t^2+1}|}}{\left(\sqrt{y^2-1}+t\right)}dt\right|^2$$
$$\leq C\int_{1<|x|<A}dx\int_{1<|y|<2}dy\left|\int_0^{A+B}
\frac{\log{|x-1|}+\log{(\sqrt{y^2-1})}}{\left(\sqrt{y^2-1}+t\right)}dt\right|^2$$
$$\leq C\int_{1<|x|<A}dx\int_{1<|y|<2}dy\log^2{\sqrt{y^2-1}}\cdot\left(\log{|x-1|}
+\log{\sqrt{y^2-1}}\right)^2<C<\infty.$$
\indent
For $t\in \left[0,A+B\right]\setminus S(x,y)$ we have
$$1+\frac{t^2}{2}\geq \sqrt{t^2+1}> x-\frac{1}{2}\left(x-1\right)\sqrt{y^2-1},$$
and therefore
$$\frac{t^2}{2}\geq \left(x-1\right)\left[1-\frac{1}{2}\sqrt{y^2-1}\right],$$
or
$$t\geq C\sqrt{x-1}.$$
\indent
Using the last inequality we obtain
$$\left|dt\right|\leq\left|\frac{\sqrt{t^2+1}}{t}du\right|
\leq C\left|\frac{1}{\sqrt{x-1}}du\right|,$$
and switching to variable $u=\sqrt{t^2+1}$ for $\left[0,A+B\right]\setminus S(x,y)$,
we obtain
$$\int_{1<|x|<A}dx\int_{1<|y|<2}dy
\left|\int_{\left\{\left[0,A+B\right]\setminus S(x,y)\right\}}
\frac{\log{|x-\sqrt{t^2+1}|}}{\left(\sqrt{y^2-1}+t\right)}dt\right|^2$$
$$\leq C\int_{1<|x|<A}dx\int_{1<|y|<2}dy
\left|\int_{x-\frac{1}{2}(x-1)\sqrt{y^2-1}}^{x+\frac{1}{2}(x-1)\sqrt{y^2-1}}\
\frac{\log{|u-x|}}{\sqrt{y^2-1}\cdot\sqrt{x-1}}du\right|^2$$
$$\leq C\int_{1<|x|<A}dx\int_{1<|y|<2}dy
\left|\frac{(x-1)\sqrt{y^2-1}\cdot\left(\log{(x-1)}+\log{\sqrt{y^2-1}}\right)}
{\sqrt{y^2-1}\cdot\sqrt{x-1}}\right|^2<C<\infty.$$
\indent
Combining the last two estimates above we obtain
$$\int_{1<|x|<A}dx\int_{1<|y|<2}dy\left|\int_0^{A+B}
\frac{\log{|x-\sqrt{t^2+1}|}}{\left(\sqrt{y^2-1}+t\right)}dt\right|^2<C<\infty,$$
and therefore
\begin{equation}\label{N2Estimate}
\int_{1<|x|<A}dx\int_{1<|y|<2}\frac{dy}{|y|^2}\left|\int_0^{A+B}M(x-\sqrt{t^2+1},\lambda)
\frac{\left(y+\sqrt{t^2+1}\right)dt}{\left(\sqrt{y^2-1}+t\right)}\right|^2
<C|\lambda\log{\lambda}|^2.
\end{equation}
\indent
To prove estimate (\ref{NEstimate}) we use the following lemma.\\

\begin{lemma}\label{MEstimates}\ The following estimates hold for $1<|x|<A$
and $\mbox{Re}\lambda\in [\sigma_1,\sigma_2]$
\begin{equation}\label{MNorm}
\begin{array}{ll}
{\dis \int_{\R}|M(x-u,\lambda)|^2u^2du<C|\lambda|^{1+\epsilon}\ \mbox{for arbitrary}\
\epsilon>0, }\vspace{0.1in}\\
{\dis \left|\int_0^{\infty}M(x-\sqrt{t^2+1},\lambda)dt\right|<C \sqrt{|\lambda|}. }
\end{array}
\end{equation}
\end{lemma}
\indent
{\bf Proof.}\ Using representation (\ref{RRepresentation})
for $1<|x|<A$ and $|\lambda(x-u)|\leq B$ we obtain
$$M(x-u,\lambda)\cdot u=\left[\lambda\log\left(\lambda|x-u|\right) \alpha(\lambda|x-u|)
+\lambda \beta\left(\lambda|x-u|\right)+\gamma\left(\lambda|x-u|\right)\right]u,$$
and therefore
$$\int_{|x-u|\leq B/|\lambda|}|M(x-u,\lambda)|^2u^2du$$
$$<C\int_{|x-u|\leq B/|\lambda|}\left[|\lambda|^2\left(|\log{\lambda}|^2+\log^2|x-u|\right)
|\alpha(\lambda|x-u|)|^2\right.$$
$$\left.+|\lambda|^2 |\beta\left(\lambda|x-u|\right)|^2
+|\gamma\left(\lambda|x-u|\right)|^2\right]u^2du
<C|\lambda||\log{\lambda}|^2.$$
\indent
For $1<|x|<A$ and $|\lambda(x-u)|\geq B$ from (\ref{RRepresentation}) we have
$$M(x-u,\lambda)=\lambda \delta\left(\lambda|x-u|\right)\frac{e^{-(\sigma+i\eta)|x-u|}}
{\sqrt{|\lambda|\cdot|x-u|}},$$
and therefore
$$\int_{|x-u|\geq B/|\lambda|}|M(x-u,\lambda)|^2u^2dx$$
$$<C\int_{\R}|\lambda|^{1+\epsilon}|\delta\left(\lambda|x-u|\right)|^2
\frac{e^{-2\sigma|x-u|}u^2du}{\left(|\lambda(x-u)|\right)^{\epsilon}|x-u|^{1-\epsilon}}
<C|\lambda|^{1+\epsilon}$$
for any $\epsilon>0$.\\
\indent
Combining the estimates above we obtain the first estimate from (\ref{MNorm}).\\
\indent
For the second integral in (\ref{MNorm}) we use representation (\ref{RRepresentation}),
and obtain for $1<|x|<A$ and $|\lambda(x-\sqrt{t^2+1})|\leq B$
$$\left|\int_{|x-\sqrt{t^2+1}|\leq B/|\lambda|}M(x-\sqrt{t^2+1},\lambda)dt\right|$$
$$\leq \int_{|x-\sqrt{t^2+1}|\leq B/|\lambda|}
\Big|\lambda\log\left(\lambda|x-\sqrt{t^2+1}|\right)\alpha(\lambda|x-\sqrt{t^2+1}|)$$
$$+\lambda \beta\left(\lambda|x-\sqrt{t^2+1}|\right)
+\gamma\left(\lambda|x-\sqrt{t^2+1}|\right)\Big|dt\leq C|\log{\lambda}|,$$
where we used the fact that the length of the interval of integration
is bounded by ${\dis \frac{C}{|\lambda|}}$ for some $C>0$.\\
\indent
For $1<|x|<A$ and $|\lambda(x-\sqrt{t^2+1})|> B$ using representation (\ref{RRepresentation})
we obtain
$$\left|\int_{|x-\sqrt{t^2+1}|> B/|\lambda|}M(x-\sqrt{t^2+1},\lambda)dt\right|$$
$$< C\int_{\R}\lambda \delta\left(\lambda|x-\sqrt{t^2+1}|\right)
\frac{e^{-(\sigma+i\eta)|x-\sqrt{t^2+1}|}dt}{\sqrt{|\lambda|\cdot|x-\sqrt{t^2+1}|}}
<C \sqrt{|\lambda|}.$$
\indent
Combining the two estimates above we obtain the second estimate of (\ref{MNorm}).
\qed

\indent
Using now estimates (\ref{MNorm}) from the lemma above in estimates (\ref{N1Estimate}),
(\ref{N2Estimate1}), and (\ref{N2Estimate2}) and combining them with estimate (\ref{N2Estimate})
we obtain estimate (\ref{NEstimate}) of Proposition~\ref{RFredholm}.
\qed

\indent
Proposition~\ref{RFredholm} allows us to reduce the question of solvability of
(\ref{MainEquation}) to the solvability of corresponding equation for ${\cal G}_{\lambda}$.
Namely, calling those $\lambda$ for which operator ${\cal G}_{\lambda}$ is not invertible
by {\it characteristic values of} ${\cal G}_{\lambda}$, we have

\begin{proposition}\label{RSolvability}\ If $\lambda_0$ is not a characteristic value
of ${\cal G}_{\lambda}$, then for arbitrary function $f \in {\cal L}^2\left(I^c(1)\right)$
and $\lambda=\lambda_0$ there exists a solution $h$ of equation (\ref{MainEquation}) such
that $h \in {\cal L}^p\left(I^c(1)\right)$ for any ${\dis p<\frac{4}{3} }$.
\end{proposition}
\indent
{\bf Proof.}\ Considering  a solution of
$${\cal G}_{\lambda}[g]={\cal R}_{\lambda}\circ{\cal P}^{-1}[g]=f$$
we define $h={\cal P}^{-1}[g]$, which satisfies equation (\ref{MainEquation})
and belongs to ${\cal L}^p\left(I^c(1)\right)$ for any
${\dis p<\frac{4}{3} }$ according to Proposition~\ref{Cofinite}.

\qed

\section{ The resolvent of operator ${\cal G}_{\lambda}$.}
\label{GLambdaResolvent}

\indent
In this section we construct the resolvent of the operator ${\cal G}_{\lambda}$
and show that it is a Fredholm operator also analytically depending on
$\lambda \in \left\{\mbox{Re}\lambda>\sigma_1\right\}$.\\
\indent
Let ${\cal T}:L^2(\R)\to L^2(\R)$ be an integral operator with kernel $T(x,y)$ satisfying
Hilbert-Schmidt condition. Following \cite{C}, we consider for operator ${\cal T}$
Hilbert's modification of the original Fredholm's determinants:
$${\cal D}_{T,m}\left(t_1,\dots,t_m\right)
=\left|
\begin{array}{cccc}
0&T(t_1,t_2)&\cdots&T(t_1,t_m)\\
T(t_2,t_1)&0&\cdots&T(t_2,t_m)\\
\vdots&&&\vdots\\
T(t_m,t_1)&\cdots&T(t_m,t_{m-1})&0
\end{array}
\right|,$$
\begin{equation}\label{D_T}
{\cal D}_{T}=1+\sum_{m=1}^{\infty}\delta_m
=1+\sum_{m=1}^{\infty}\frac{1}{m!}\int_{\R}\cdots\int_{\R}
{\cal D}_{T,m}\left(t_1,\dots,t_m\right)dt_1\cdots dt_m,
\end{equation}
$${\cal D}_{T,m}\left(\begin{array}{c}
x\\
y
\end{array}t_1,\dots,t_m\right)
=\left|
\begin{array}{cccc}
T(x,y)&T(x,t_1)&\cdots&T(x,t_m)\\
T(t_1,y)&0&\cdots&T(t_1,t_m)\\
\vdots&&&\vdots\\
T(t_m,y)&\cdots&T(t_m,t_{m-1})&0
\end{array}
\right|,$$
and
\begin{equation}\label{D_Tfunction}
\begin{array}{ll}
{\dis {\cal D}_{T}\left(\begin{array}{c}
x\\
y
\end{array}\right)
=T(x,y)+\sum_{m=1}^{\infty}\delta_m\left(\begin{array}{c}
x\\
y
\end{array}\right) }\vspace{0.1in}\\
{\dis =T(x,y)+\sum_{m=1}^{\infty}\frac{1}{m!}\int_{\R}\cdots\int_{\R}
{\cal D}_{T,m}\left(
\begin{array}{c}
x\\
y
\end{array}t_1,\dots,t_m\right)dt_1\cdots dt_m. }
\end{array}
\end{equation}
\indent
We start with the following proposition, which summarizes the results from \cite{C}
(cf. also \cite{M}), that will be used in the construction of the resolvent of
${\cal G}_{\lambda}$.

\begin{proposition}\label{Resolvent} (\cite{C})\ Let function
$T(x,y):\R^2\to \C$ satisfy Hilbert-Schmidt condition
$$\|T\|^2=\int_{\R^2}\left|T(x,y)\right|^2dxdy<\infty.$$
\indent
Then function
${\cal D}_{T}\left(\begin{array}{c}
x\\
y
\end{array}\right)\in L^2(\R^2)$ is well defined, and the following estimates hold:
\begin{equation}\label{DEstimate1}
\left|\delta_m\right|\leq\left(\frac{e}{m}\right)^{m/2}\|T\|^m,\
\left|{\cal D}_{T}\right|\leq e^{\frac{\|T\|^2}{2}},
\end{equation}
\begin{equation}\label{DEstimate2}
\left|{\cal D}_{T}\left(\begin{array}{c}
x\\
y
\end{array}\right)\right|
\leq e^{\frac{\|T\|^2}{2}}
\left(|T(x,y)|+\sqrt{e}\alpha(x)\beta(y)\right),
\end{equation}
where
$$\alpha^2(x)=\int_{\R}\left|T(x,t)\right|^2dt,\hspace{0.1in}
\beta^2(y)=\int_{\R}\left|T(t,y)\right|^2dt.$$
\indent
If ${\cal D}_{T}\neq 0$
then kernel
\begin{equation}\label{HResolvent}
H(x,y)=\left[{\cal D}_{T}\right]^{-1}
\cdot{\cal D}_{T}\left(\begin{array}{c}
x\\
y
\end{array}\right)
\end{equation}
defines the resolvent of operator ${\cal I}-{\cal T}$, i.e. it satisfies the
following equations
\begin{equation}\label{ResolventEquations}
\begin{array}{ll}
{\dis H(x,y)+\int_{\R}T(x,t)\cdot H(t,y)dt=T(x,y), }
\vspace{0.1in}\\
{\dis H(x,y)+\int_{\R}T(t,y)\cdot H(x,t)dt=T(x,y), }
\end{array}
\end{equation}
and therefore operator ${\cal I}-{\cal H}$ is the inverse of
operator ${\cal I}+{\cal T}$.
\end{proposition}
\qed

\indent
Using Proposition~\ref{Resolvent}, we construct the resolvent of
operator ${\cal G}_{\lambda}=\pi U\left({\cal I}+{\cal N}_{\lambda}\right)$,
defined in (\ref{GOperator}), and prove the estimate that will be necessary in the proof
of Theorem~\ref{Main}.\\

\begin{proposition}\label{GResolvent}\ The set of characteristic values of operator
${\cal G}_{\lambda}$ coincides with the set
$$E({\cal G})=\left\{\lambda\in \C:\ \mbox{Re}\lambda>\sigma_1,\ {\cal D}_{N_{\lambda}}=0\right\}$$
and consists of at most countably many isolated points.\\
\indent
For $\lambda \notin E({\cal G})$ there exists an operator ${\cal H}_{\lambda}$ with kernel
$H(x,y,\lambda)$ satisfying the Hilbert-Schmidt condition and such that
operator ${\cal I}-{\cal H}_{\lambda}$ is the inverse of operator ${\cal I}+{\cal N}_{\lambda}$,
and therefore operator ${\dis \frac{1}{\pi U}\left({\cal I}-{\cal H}_{\lambda}\right) }$
is the inverse of operator ${\cal G}_{\lambda}$.\\
\indent
If function ${\cal D}_{N}(\lambda)={\cal D}_{N_{\lambda}}$ has no zeros in a strip
$\left\{\lambda:\ \sigma_1<\mbox{Re}\lambda<\sigma_2\right\}$, then operator ${\cal H}_{\lambda}$
admits estimate
\begin{equation}\label{HNormEstimate}
\left||{\cal H}_{\lambda}\right\|
<\exp\left\{e^{|\eta|}\cdot(1+|\eta|)^{2+\epsilon}\right\}
\end{equation}
for $\lambda\in \left\{\sigma_1+\gamma <\mbox{Re}\lambda <\sigma_2-\gamma\right\}$
and arbitrary $\epsilon>0$.
\end{proposition}
\indent
{\bf Proof.}\ Applying Proposition~\ref{Resolvent} to operator ${\cal N}_{\lambda}$
we obtain the existence of functions
\begin{equation}\label{D_N}
{\cal D}_{N}(\lambda)={\cal D}_{N_{\lambda}}
\end{equation}
and
$${\cal D}_{N}\left(\left.\begin{array}{c}
x\\
y
\end{array}\right|\lambda\right)
={\cal D}_{N_{\lambda}}\left(\begin{array}{c}
x\\
y
\end{array}\right)$$
such that for any fixed $\lambda$, satisfying ${\cal D}_{N}(\lambda)\neq 0$, kernel
$$H(x,y,\lambda)=\left[{\cal D}_{N}(\lambda)\right]^{-1}
\cdot{\cal D}_{N}\left(\left.\begin{array}{c}
x\\
y
\end{array}\right|\lambda\right)\in L^2(\R^2),$$
and operator ${\cal I}-{\cal H}_{\lambda}$ is the inverse of
operator ${\cal I}+{\cal N}_{\lambda}$.\\
\indent
Terms of the series (\ref{D_T}) for ${\cal N}_{\lambda}$ analytically depend on $\lambda$,
and according to estimates (\ref{DEstimate1}) this series converges uniformly with respect
to $\lambda$ on compact subsets of $\left\{\lambda\in\C:\ \mbox{Re}\lambda>\sigma_1\right\}$.
Therefore, ${\cal D}_{N}(\lambda)$ is an analytic function on
$\left\{\lambda\in\C:\ \mbox{Re}\lambda>\sigma_1\right\}$, and the set $E({\cal G})$
consists of at most countably many isolated points.\\
\indent
Analyticity of ${\cal I}-{\cal H}_{\lambda}$ with respect to $\lambda$ on
$$\left\{\lambda\in\C:\ \mbox{Re}\lambda>\sigma_1\right\}\setminus E({\cal G})$$
follows from the Theorem VI.14 in \cite{RS}. It is proved by approximation of the kernel by
degenerate kernels and by the argument that can be traced back to at least \cite{M}.\\
\indent
To prove estimate (\ref{HNormEstimate}) we use the well known estimate (\cite{L})
$$\left\|T\right\|^2\leq \int_{\R}|T(x,y)|^2dxdy$$
for integral operators. Using this estimate, estimates (\ref{DEstimate2}) and
(\ref{NEstimate}) we obtain
$$\left\|{\cal D}_{N}\left(\left.\begin{array}{c}
x\\
y
\end{array}\right|\lambda\right)\right\|
< \exp{\left\{C(1+|\eta|)^2\cdot\log^2{|\eta|}\right\}}
(1+|\eta|)^4\cdot\log^4{|\eta|}.$$
\indent
To estimate function $\left[{\cal D}_{N}(\lambda)\right]^{-1}$ for
$\lambda\in \left\{\sigma_1+\gamma <\mbox{Re}\lambda <\sigma_2-\gamma\right\}$
we use the following lemma.

\begin{lemma}\label{Estimate-ofD-1}\ If function ${\cal D}_{N}(\lambda)={\cal D}_{N_{\lambda}}$
has no zeros in the strip
$\left\{\lambda:\ \sigma_1<\mbox{Re}\lambda<\sigma_2\right\}$, then estimate
\begin{equation}\label{D-1Estimate}
\left|1/{\cal D}_{N}\left(\lambda\right)\right|
<\exp\left\{e^{|\eta|}\cdot(1+|\eta|)^{2+\epsilon}\right\}
\end{equation}
holds for $\lambda\in \left\{\sigma_1+\gamma <\mbox{Re}\lambda <\sigma_2-\gamma\right\}$
with fixed $\gamma>0$ and arbitrary $\epsilon>0$.
\end{lemma}
\indent
{\bf Proof.}\ We consider a biholomorphic map
$$\Psi:\ \left\{\lambda:\ \sigma_1<\mbox{Re}\lambda<\sigma_2\right\}\to
\D(1)=\left\{z\in\C: |z|<1\right\},$$
defined by the formula
$$\Psi(\lambda)=\frac{e^{i(\lambda-\sigma_1)\frac{\pi}{\sigma_2-\sigma_1}}-i}
{e^{i(\lambda-\sigma_1)\frac{\pi}{\sigma_2-\sigma_1}}+i}.$$
\indent
Denoting
$$w=u+iv=e^{i(\lambda-\sigma_1)\frac{\pi}{\sigma_2-\sigma_1}},$$
we obtain for the circle $C(r)=\left\{z:\ |z|=r\right\}$
$$\Psi^{-1}\left(C(r)\right)=\left\{\sigma+i\eta:\
\left|e^{i(\lambda-\sigma_1)\frac{\pi}{\sigma_2-\sigma_1}}-i\right|
=r\left|e^{i(\lambda-\sigma_1)\frac{\pi}{\sigma_2-\sigma_1}}+i\right|\right\}$$
$$=\left\{u+iv:\ \left(u^2+v^2-2v+1\right)
=r^2\left(u^2+v^2+2v+1\right)\right\}$$
$$=\left\{u+iv:\ u^2+\left(v-\frac{1+r^2}{1-r^2}\right)^2
=\frac{4r^2}{(1-r^2)^2}\right\}.$$
\indent
Introducing coordinates
$$t=\mbox{Re}\frac{\pi(\lambda-\sigma_1)}{\sigma_2-\sigma_1},\
s=\mbox{Im}\frac{\pi(\lambda-\sigma_1)}{\sigma_2-\sigma_1},$$
such that
$$w=u+iv=e^{i(\lambda-\sigma_1)\frac{\pi}{\sigma_2-\sigma_1}}
=e^{it-s}=e^{-s}\left(\cos{t}+\sin{t}\right),$$
we can rewrite the last condition as a quadratic equation with respect
to $e^{-s}$ for fixed $t$
$$\left(e^{-s}-\sin{t}\frac{1+r^2}{1-r^2}\right)^2
+\cos^2{t}\left(\frac{1+r^2}{1-r^2}\right)^2-\frac{4r^2}{(1-r^2)^2}=0.$$
\indent
Solving equation above we obtain
$$e^{-s}=\sin{t}\frac{1+r^2}{1-r^2}
\pm\sqrt{\frac{4r^2}{(1-r^2)^2}-\cos^2{t}\left(\frac{1+r^2}{1-r^2}\right)^2}$$
with solutions existing for $t$ such that
$$|\cos{t}|\leq \frac{2r}{1-r^2}\frac{1-r^2}{1+r^2}=\frac{2r}{1+r^2}.$$
\indent
The maximal value for $e^{-s}$ is achieved at $t=\frac{\pi}{2}$ and it is
$$e^{-s}=\frac{1+r^2}{1-r^2}+\frac{2r}{1-r^2}
=\frac{1+r^2+2r}{1-r^2}=\frac{(1+r)^2}{1-r^2}=\frac{1+r}{1-r}.$$
Therefore the maximal value for $|s|$ is achieved at $t=\frac{\pi}{2}$, is
equal to $|s|=\log\left(\frac{1+r}{1-r}\right)$, and for $r=1-\delta$
we have the maximal value
\begin{equation}\label{s_max}
\max|s|=\log\left(\frac{1+r}{1-r}\right)=-\log{\delta}+\log{(2-\delta)}.
\end{equation}
\indent
Since function ${\cal D}_{N}\left(\lambda\right)$ has no zeros
in $\left\{\lambda:\ \sigma_1<\mbox{Re}\lambda<\sigma_2\right\}$ we can consider
analytic function $\log\left({\cal D}_{N}\left(\lambda\right)\right)$
in this strip, and using estimates (\ref{DEstimate1}) and (\ref{NEstimate}),
and equality (\ref{s_max}), we obtain the following estimate for
$z=(1-\delta)e^{i\theta}$
$$\log\left|{\cal D}_{N}\left(\Psi^{-1}(z)\right)\right|
\leq \frac{\left\|N_{\Psi^{-1}(z)}\right\|^2}{2}
\leq C\left|\Psi^{-1}(z)\cdot\log\left(\Psi^{-1}(z)\right)\right|^2$$
$$\leq C|\log{\delta}\cdot\log\left(\log{\delta}\right)|^2.$$
\indent
Using then the Borel-Caratheodory inequality (\cite{Ti1}, \cite{Boa}) on disks with radii
$$1-2\delta=r<R=1-\delta,$$
we obtain
$$\left|\log\left({\cal D}_{N}\left(\Psi^{-1}(z)\right)\right)
\right|_{\{|z|=1-2\delta\}}$$
$$\leq\frac{2-4\delta}{\delta}\max_{|z|=R}
{\mbox{Re}\left\{\log\left({\cal D}_{N}\left(\Psi^{-1}(z)\right)\right)\right\}}
+\frac{1-\delta+1-2\delta}{\delta}|\log\left({\cal D}_{N}\left(\Psi^{-1}(0)\right)\right)|$$
$$< \frac{C}{\delta}\log^2{\delta}\cdot\log^2\left(\log{\delta}\right),$$
or
$$-\frac{C}{\delta}\log^2{\delta}\cdot\log^2\left(\log{\delta}\right)
<\mbox{Re}\left\{\log\left({\cal D}_{N}\left(\Psi^{-1}(z)\right)\right)
\right\}\Big|_{\{|z|=1-2\delta\}}
<\frac{C}{\delta}\log^2{\delta}\cdot\log^2\left(\log{\delta}\right).$$
\indent
From the last estimate we obtain an estimate for
the function $\left|1/{\cal D}_{N}\left(\Psi^{-1}(z)\right)\right|$ in the disk
$\D(1-2\delta)$:
\begin{equation}\label{D-lowestimate}
\left|1/{\cal D}_{N}\left(\Psi^{-1}(z)\right)\right|\Bigg|_{\{|z|\leq 1-2\delta\}}
< \exp\left\{\frac{|\log{\delta}|^{2+\epsilon}}{\delta}\right\}
\end{equation}
for arbitrary $\epsilon>0$.\\
\indent
For a fixed $t\in\left(0,\pi\right)$ and arbitrary $s$ we have that
$t+is\in\Psi^{-1}\left(\D(r)\right)$ with $r=1-2\delta$ if
$$e^{|s|}\leq \sin{t}\cdot\frac{1+r^2}{1-r^2}
+\sqrt{\frac{4r^2}{(1-r^2)^2}-\cos^2{t}\cdot\left(\frac{1+r^2}{1-r^2}\right)^2}$$
$$=\sin{t}\cdot\frac{2-4\delta+4\delta^2}{2\delta(2-2\delta)}
+\frac{\sqrt{4(1-2\delta)^2-\cos^2{t}\cdot(2-4\delta+4\delta^2)^2}}
{2\delta(2-2\delta)},$$
and therefore for any interval $\left[\gamma^{\prime},\pi-\gamma^{\prime}\right]$
there exist constants $C_1,\ C_2$ such that conditions
$$t\in \left[\gamma^{\prime},\pi-\gamma^{\prime}\right],\hspace{0.1in}
\frac{C_1}{\delta}<e^{|s|}<\frac{C_2}{\delta}$$
imply that $t+is \in \Psi^{-1}\left(\D(1-2\delta)\right)$.\\
\indent
Using then estimate (\ref{D-lowestimate}) we obtain for $\lambda$ with
$\mbox{Re}\lambda \in \left[\sigma_1+\frac{\gamma^{\prime}(\sigma_2-\sigma_1)}{\pi},
\sigma_2-\frac{(\pi-\gamma^{\prime})(\sigma_2-\sigma_1)}{\pi}\right]$ the estimate
$$\left|1/{\cal D}_{N}\left(\lambda\right)\right|
<\exp\left\{e^{|s|}\cdot(1+|s|)^{2+\epsilon}\right\}$$
for arbitrary $\epsilon>0$, which leads to estimate (\ref{D-1Estimate}).
\qed

\indent
Combining now estimate for ${\dis \left\|{\cal D}_{N}\left(\left.\begin{array}{c}
x\\
y
\end{array}\right|\lambda\right)\right\| }$ with (\ref{D-1Estimate}) we obtain
estimate (\ref{HNormEstimate}).
\qed

\section{Proof of Theorem~\ref{Main}.}
\label{Proof}

\indent
Before proving Theorem~\ref{Main} we will prove two lemmas, that will be used in
the proof of this theorem.\\
\indent
In order to assure applicability of Proposition~\ref{RSolvability} to
$f_a$, defined in (\ref{fDefinition}), we have to prove that
$$f_a \in {\cal L}^2\left(I^c(1)\right)$$
for $\widehat{w_a}$ satisfying (\ref{wCondition}). In the lemma below we
prove the necessary property of $f_a$.

\begin{lemma}\label{fCondition}\ If $\widehat{w_a}$ satisfies condition
(\ref{wCondition}) then $f_a(x,\lambda)$ defined by the formula (\ref{fDefinition})
is a function in ${\cal L}^2\left(I^c(1)\right)$ for any fixed $\lambda$,
which satisfies the estimate
\begin{equation}\label{f_aEstimate}
\left\|f_a(\cdot,\sigma+i\eta)\right\|_{{\cal L}^2\left(I^c(1)\right)}
<C\exp{\left\{-e^{|\eta|}\cdot(1+|\eta|)^{2+\epsilon}\right\}}
\end{equation}
with some $\epsilon>0$ for $\sigma\in[\sigma_1,\sigma_2]$.
\end{lemma}
\indent
{\bf Proof.}\ For a fixed $\lambda=\sigma+i\eta$ with $\sigma\in[\sigma_1,\sigma_2]$ we choose
$B>1$, and using second representation from (\ref{RRepresentation}) of $R(x,\lambda)$
for $|\lambda x|>B$, obtain an estimate
$$\left|R(x-y,\lambda)\right|<C\frac{|\lambda|^{1/2}e^{-\lambda|x-y|}}{\sqrt{|x-y|}}.$$
Using then condition (\ref{wCondition}), we have
\begin{equation}\label{BigxEstimate}
\left(\int_{|x|>B/|\lambda|}\left|f_a(x,\lambda)\right|^2dx\right)^{1/2}
\end{equation}
$$=\frac{1}{\pi^2}\left|\int_{|x|>B/|\lambda|}\left|\int_{-1}^1 e^{-d(\lambda)y}
R(x-y,\lambda)\widehat{w_a}(y,\lambda)dy\right|^2dx\right|^{1/2}$$
$$< C|\lambda|^{1/2}\left(\int_{|x|>B/|\lambda|}\left(\int_{-1}^1
e^{-\sigma|x|}\left|\widehat{w_a}(y,\lambda)\right|dy\right)^2dx\right)^{1/2}$$
$$< C|\lambda|^{1/2}\int_{-1}^1\left|\widehat{w_a}(y,\lambda)\right|dy
<C|\lambda|^{1/2}\int_{-1}^1\left|\widehat{w_a}(y,\lambda)\right|^2dy$$
$$<C\exp{\left\{-e^{|\eta|}\cdot(1+|\eta|)^{2+\epsilon}\right\}}.$$
\indent
For $|\lambda x|<B$ we use the first representation from (\ref{RRepresentation})
for $R(x-y,\lambda)$. Since the Hilbert transform is a bounded linear operator from
$L^q$ into $L^q$ (see \cite{Ti2}, \cite{Tr}), and kernels $\alpha\left(\lambda(x-y)\right)$,
$\beta\left(\lambda(x-y)\right)$, and  $\gamma\left(\lambda(x-y)\right)$ from
(\ref{RRepresentation}) are bounded, we obtain
\begin{equation}\label{SmallxEstimate}
\left(\int_{|x|<|B/\lambda|}\left|f_a(x,\lambda)\right|^2dx\right)^{1/2}
<C\left|\lambda\log{\lambda}\right|
\cdot\left\|\widehat{w_a}(y,\lambda)\right\|_{L^2\left(I(1)\right)}
\end{equation}
$$<C\exp{\left\{-e^{|\eta|}\cdot(1+|\eta|)^{2+\epsilon}\right\}},$$
where in the last inequality we used condition (\ref{wCondition}).\\
\indent
Combining estimates (\ref{BigxEstimate}) and (\ref{SmallxEstimate}) we obtain
(\ref{f_aEstimate}).
\qed

\begin{lemma}\label{LaplaceInverse}\ If a function $h(y,\lambda)$ satisfies
estimate
\begin{equation}\label{hCondition}
\int_{-\infty}^{\infty}e^{-\sigma_1\cdot|y|}|h(y,\sigma+i\eta)|dy
<\frac{C}{(1+|\eta|)^{\frac{5}{2}+\epsilon}}
\end{equation}
for some $\epsilon>0$ and $\sigma_1<\mbox{Re}\lambda<\sigma_2$, then function
\begin{equation}\label{xi_inL}
\xi(x,z,\lambda)=e^{d(\sigma+i\eta)x}\left[\int_{-\infty}^\infty
K_0\left(r(\sigma+i\eta)\left(\frac{(x-y)^2}{1-M^2}+z^2\right)^{\frac{1}{2}}\right)
h(y,\sigma+i\eta)dy\right]\in{L^1_{\eta}(\R)}
\end{equation}
for $\sigma\in [\sigma_1,\sigma_2]$, and satisfies conditions
\begin{equation}\label{xiDerivatives}
\frac{\partial^2\xi(x,z,\sigma+i\eta)}{\partial x^2},\
\frac{\partial^2\xi(x,z,\sigma+i\eta)}{\partial z^2},\
|\eta|^2\xi(x,z,\sigma+i\eta),\
|\eta|\frac{\partial\xi(x,z,\sigma+i\eta)}{\partial x}\in L^1(\R_{\eta}).
\end{equation}
\indent
Function
$$\phi(x,z,t)=-\frac{1}{2\pi\sqrt{1-M^2}}
\int_{\sigma-i\infty}^{\sigma+i\infty}e^{d(\sigma+i\eta)x}$$
$$\times\left[\int_{-\infty}^\infty K_0\left(r(\sigma+i\eta)
\left(\frac{(x-y)^2}{1-M^2}+z^2\right)^{\frac{1}{2}}\right)
h(y,\sigma+i\eta)dy\right]e^{(\sigma+i\eta)t}d\eta$$
is then well defined for $z>0$, and doesn't depend on $\sigma\in [\sigma_1,\sigma_2]$.
\end{lemma}
\indent
{\bf Proof.}\ To prove inclusion (\ref{xi_inL}) of the Lemma it suffices to prove that
under conditions above estimate
\begin{equation}\label{LaplaceEstimate1}
\left\|e^{d(\sigma+i\eta)x}\left[\int_{-\infty}^\infty K_0\left(r(\sigma+i\eta)
\left(\frac{(x-y)^2}{1-M^2}+z^2\right)^{\frac{1}{2}}\right)
h(y,\sigma+i\eta)dy\right]\right\|_{L^1_{\eta}(\R)}<C(M,z)
\end{equation}
holds uniformly with respect to $\sigma \in \left[\sigma_1,\sigma_2\right]$
for fixed $x$, fixed $z>0$, and for some $\sigma_1>\sigma_a$. Applying then Theorem 47 from
\cite{Boc} we will obtain the second part of the Lemma.\\
\indent
Using asymptotics of $K_0(\zeta)$ for large and for small $|\zeta|$ (\cite{EMOT}) we obtain
the existence for fixed $z>0$ of a constant $A(z)>0$, large enough, such that estimates
\begin{equation}\label{LargeAEstimate}
\left|K_0\left(r(\sigma+i\eta)
\left(\frac{(x-y)^2}{1-M^2}+z^2\right)^{\frac{1}{2}}\right)\right|
<\frac{C(M,z)e^{-\sigma|x-y|}}{\sqrt{|\sigma+i\eta|\cdot|x-y|}}\ \mbox{for}\ |x-y|>A(z),
\end{equation}
and
\begin{equation}\label{SmallAEstimate}
\left|K_0\left(r(\sigma+i\eta)
\left(\frac{(x-y)^2}{1-M^2}+z^2\right)^{\frac{1}{2}}\right)\right|
<\frac{C(M,z)}{\sqrt{|\sigma+i\eta|}}\ \mbox{for}\ |x-y|<A(z),
\end{equation}
hold uniformly for $\sigma \in \left[\sigma_1,\sigma_2\right]$,
with a constant $C$ depending on $M$ and $z$.\\
\indent
Combining estimates (\ref{LargeAEstimate}) and (\ref{SmallAEstimate}) with the
estimate for $h(y,\lambda)$ we obtain
$$\left|e^{d(\sigma+i\eta)x}\left[\int_{-\infty}^\infty K_0\left(r(\sigma+i\eta)
\left(\frac{(x-y)^2}{1-M^2}+z^2\right)^{\frac{1}{2}}\right)
h(y,\sigma+i\eta)dy\right]\right|$$
$$<\frac{C(M,z)}{\sqrt{|\sigma+i\eta|}}
\int_{-\infty}^\infty e^{-\sigma_1\cdot|y|}|h(y,\sigma+i\eta)|dy
<\frac{C(M,z)}{\left(1+|\eta|\right)^{3+\epsilon}},$$
for $z>0$, which leads to estimate (\ref{LaplaceEstimate1}).\\
\indent
Again using estimates (\ref{LargeAEstimate}) and (\ref{SmallAEstimate}) and analogous
estimates for
$$\frac{\partial}{\partial x}K_0\left(r(\sigma+i\eta)
\left(\frac{(x-y)^2}{1-M^2}+z^2\right)^{\frac{1}{2}}\right),\
\frac{\partial^2}{\partial x^2}K_0\left(r(\sigma+i\eta)
\left(\frac{(x-y)^2}{1-M^2}+z^2\right)^{\frac{1}{2}}\right),$$
and
$$\frac{\partial^2}{\partial z^2}K_0\left(r(\sigma+i\eta)
\left(\frac{(x-y)^2}{1-M^2}+z^2\right)^{\frac{1}{2}}\right)$$
we obtain inclusions (\ref{xiDerivatives}).
\qed\\

\indent
To prove Theorem~\ref{Main} we consider $w_a$ satisfying condition (\ref{wCondition}),
and define $f_a$ by the formula (\ref{fDefinition}). Using Lemma~\ref{fCondition} we obtain
that $f_a$ satisfies estimate (\ref{f_aEstimate}). Applying Proposition~\ref{RSolvability}
to $f_a$ and using estimate (\ref{HNormEstimate}) from Proposition~\ref{GResolvent}
we obtain the existence of $h_a$ satisfying equation (\ref{MainEquation}) and such that
$$\left\|h_a(\cdot,\sigma+i\eta)\right\|_{{\cal L}^p\left(I^c(1)\right)}
<\exp\left\{e^{|\eta|}\cdot(1+|\eta|)^{2+\epsilon}\right\}
\cdot\left\|f_a(\cdot,\sigma+i\eta)\right\|_{{\cal L}^2\left(I^c(1)\right)}
<\frac{C}{\left(1+|\eta|\right)^m}$$
for arbitrary $m$, arbitrary $p<\frac{4}{3}$, and $\sigma\in[\sigma_1,\sigma_2]$,
with $\sigma_a<\sigma_1$.\\
\indent
Using the estimate above for $p=1$, we obtain
\begin{equation}\label{h_aFarEstimate}
\int_{|x|>1}\left|h_a(x,\sigma+i\eta)\right|\cdot|x|^{-1}dx
<\frac{C}{\left(1+|\eta|\right)^m}.
\end{equation}
\indent
From the definition of $h_a$ on $[-1,1]$ as
$$h_a(x,\lambda)=\frac{e^{-d(\lambda)x}\cdot \widehat{w_a}(x,\lambda)}{\pi}$$
and from condition (\ref{wCondition}) we obtain
$$\left\|h_a(x,\sigma+i\eta)\right\|_{L^p\left(I(1)\right)}
=\left\|e^{-d(\lambda)x}\cdot\widehat{w_a}(x,\sigma+i\eta)\right\|_{L^p\left(I(1)\right)}$$
$$<C\left\|\widehat{w_a}(\cdot,\sigma+i\eta)\right\|_{L^2\left(I(1)\right)}
<\frac{C}{\left(1+|\eta|\right)^m}\
\mbox{for}\ p<\frac{4}{3},\ \sigma\in[\sigma_1,\sigma_2]\ \mbox{with}\ \sigma_a<\sigma_1,$$
and therefore
\begin{equation}\label{h_aCloseEstimate}
\left\|h_a(\cdot,\sigma+i\eta)\right\|_{L^1\left(I(1)\right)}
<\frac{C}{\left(1+|\eta|\right)^m}
\end{equation}
for arbitrary $m>0$.\\
\indent
From the estimates (\ref{h_aFarEstimate}) and (\ref{h_aCloseEstimate}) we conclude
that function $h_a$ satisfies estimate (\ref{hCondition}), and therefore, applying
Lemma~\ref{LaplaceInverse} and Proposition~\ref{LaplaceSolution}, we obtain that function
$\phi(x,z,t)$ in formula (\ref{Solution})
is well defined and satisfies equation (\ref{LinearEquation}).\\
\indent
To prove that $\phi(x,z,t)$ satisfies boundary condition (\ref{Flow_Tangency}) we
fix $x\in[-1,1]$ and denote $\delta=\min\left\{x+1,1-x\right\}$. Then we have
\begin{equation}\label{Tangency1}
\lim_{z\to 0}\frac{\partial}{\partial z}\xi(x,z,\lambda)
=\lim_{z\to 0}\frac{\partial}{\partial z}\int_{-\infty}^\infty
S(x-y,z,\lambda)e^{d(\lambda)y}h_a(y,\lambda)dy
\end{equation}
$$=\lim_{z\to 0}\frac{\partial}{\partial z}\int_{x-\frac{\delta}{2}}^{x+\frac{\delta}{2}}
S(x-y,z,\lambda)e^{d(\lambda)y}h_a(y,\lambda)dy$$
$$+\lim_{z\to 0}\frac{\partial}{\partial z}
\int_{\R\setminus[x-\frac{\delta}{2},x+\frac{\delta}{2}]}
S(x-y,z,\lambda)e^{d(\lambda)y}h_a(y,\lambda)dy.$$
For the first integral in the right hand side of (\ref{Tangency1}) we obtain
using Lemma~\ref{Fourier}
$$\lim_{z\to 0}\frac{\partial}{\partial z}\int_{x-\frac{\delta}{2}}^{x+\frac{\delta}{2}}
S(x-y,z,\lambda)e^{d(\lambda)y}h_a(y,\lambda)dy$$
$$=-\lim_{z\to 0}\frac{\partial}{\partial z}\int_{x-\frac{\delta}{2}}^{x+\frac{\delta}{2}}
e^{d(\lambda)y}h_a(y,\lambda)dy\int_{-\infty}^\infty e^{i(x-y)\omega}
\frac{e^{\dis -z\left((1-M^2)(\omega+id(\lambda))^2+r^2(\lambda)\right)^{\frac{1}{2}}} }
{2\sqrt{(1-M^2)(\omega+id(\lambda))^2+r^2(\lambda)}}d\omega$$
$$=\widehat{w_a}(x,\lambda).$$
\indent
For the second integral in the right hand side of (\ref{Tangency1}) we have
$$\lim_{z\to 0}\frac{\partial}{\partial z}
\int_{\R\setminus[x-\frac{\delta}{2},x+\frac{\delta}{2}]}
S(x-y,z,\lambda)e^{d(\lambda)y}h_a(y,\lambda)dy$$
$$=-\frac{e^{d(\lambda)x}}{\sqrt{1-M^2}}\lim_{z\to 0}
\int_{\R\setminus[x-\frac{\delta}{2},x+\frac{\delta}{2}]}
\left[\frac{\partial}{\partial z}K_0\left(r(\lambda)
\left(\frac{(x-y)^2}{1-M^2}+z^2\right)^{\frac{1}{2}}\right)\right]h_a(y,\lambda)dy.$$
\indent
Using then estimate (\ref{h_aFarEstimate}) and equality
$$\lim_{z\to 0}\left[\frac{\partial}{\partial z}K_0\left(r(\lambda)
\left(\frac{(x-y)^2}{1-M^2}+z^2\right)^{\frac{1}{2}}\right)\right]=0$$
for $y\in \R\setminus[x-\frac{\delta}{2},x+\frac{\delta}{2}]$ we obtain that
the second integral in the right hand side of (\ref{Tangency1}) is equal to zero.\\
\indent
From equalities above we conclude that
$$\lim_{z\to 0}\frac{\partial}{\partial z}\xi(x,z,\lambda)=\widehat{w_a}(x,\lambda)$$
for $x\in [-1,1]$ and $\mbox{Re}\lambda\in [\sigma_1,\sigma_2]$, and therefore
$$\lim_{z\to 0}\frac{\partial}{\partial z}\phi(x,z,t)=w_a(x,t).$$
\indent
Straightforward substitution of $v_a(x,\lambda)=e^{d(\lambda)x}h_a(x,\lambda)$
into the formula (\ref{LapSolution}), with $h_a(x,\lambda)$ defined as
$$h_a(x,\lambda)=\left\{
\begin{array}{ll}
{\dis \frac{1}{\pi}e^{-d(\lambda)x}\cdot \widehat{w_a}(x,\lambda)\
\mbox{for}\ x\in[-1,1], }\vspace{0.1in}\\
\mbox{solution of equation}\ (\ref{MainEquation})\ \mbox{for}\ x\in\R\setminus[-1,1],
\end{array}\right.$$
shows that $\xi(x,z,\lambda)$ defined by this formula satisfies equation (\ref{ZeroCondition})
for $1<|x|<A$. Then for $\phi(x,z,t)$ defined by formula (\ref{InverseLaplace}) we will have
$$\frac{\partial\phi(x,0,t)}{\partial t}+U\frac{\partial\phi(x,0,t)}{\partial x}$$
$$=\frac{1}{2\pi}\left(\frac{\partial}{\partial t}+U\frac{\partial}{\partial x}\right)
\int_{\sigma-i\infty}^{\sigma+i\infty}e^{\lambda t}\xi(x,0,\lambda)d\eta$$
$$=\frac{1}{2\pi}\int_{\sigma-i\infty}^{\sigma+i\infty}e^{\lambda t}
\left(\lambda+U\frac{\partial}{\partial x}\right)\xi(x,0,\lambda)d\eta=0$$
for $1<|x|<A$.

\qed


\begin{thebibliography}{WIDT}
\bibitem[Ba1]{Ba1} A.V. Balakrishnan, Semigroup theory in aeroelasticity, Progress in Nonlinear
Differential Equations and Their Applications, v.42, 15-24, Birkh\"auser Verlag, 2000.
\bibitem[Ba2]{Ba2} A.V. Balakrishnan, Possio integral equation of aeroelasticity theory,
Journal of Aerospace Engineering, 16:4 (2003), 139-154.
\bibitem[BAH]{BAH} R.L. Bisplinghoff, H. Ashley, R.L. Halfman, Aeroelasticity, Dover,
New York, 1996.
\bibitem[Boa]{Boa} R. Boas, Entire functions, Academic Press, New York, 1954.
\bibitem[Boc]{Boc} S. Bochner, Lectures on Fourier integrals, Princeton University Press,
Princeton, NJ, 1959.
\bibitem[C]{C} T. Carleman, Zur Theorie der linearen Integralgleichungen, Math. Zeitschrift,
v. 9, 196-217, 1921.
\bibitem[EMOT]{EMOT} A. Erd\'elyi, W. Magnus, F. Oberhettinger, F.G. Tricomi, v.I Tables of
integral transforms, v.II Higher transcendental functions, CalTech Bateman Manuscript Project,
McGraw-Hill, 1954.
\bibitem[GR]{GR} I.S. Gradshteyn, I.M. Ryzhik, Table of integrals, series, and products,
Academic Press, 1994.
\bibitem[L]{L} P.D. Lax, Functional analysis, John Wiley \& Sons, 2002.
\bibitem[M]{M} S.G. Mikhlin, Integral equations, Pergamon Press, 1957.
\bibitem[RS]{RS} M. Reed, B. Simon, Functional Analysis, Academic Press, 1980.
\bibitem[So]{So} H. S\"ohngen, Die L\"osungen der Integralglechung und deren
Anwendung in der Tragfl\"ugeltheorie, Math. Zeitschrift, v. 45, 245-264, 1939.
\bibitem[St]{St} E. Stein, Singular integrals and differentiability properties of functions,
Princeton University Press, Princeton, NJ, 1970.
\bibitem[Ti1]{Ti1} E.C. Titchmarsh, The theory of functions, Oxford University Press, 1939.
\bibitem[Ti2]{Ti2} E.C. Titchmarsh, Introduction to the theory of Fourier integrals, Chelsea,
New York, 1986.
\bibitem[Tr]{Tr} F.G. Tricomi, Integral equations, Intersciense Publishers, New York, 1957.
\end{thebibliography}
\end{document}